\documentclass[12pt,twoside]{article}
\usepackage[english]{babel}
\usepackage[latin1]{inputenc}
\usepackage{amsmath}
\usepackage{amssymb,amsfonts}
\usepackage{graphicx}                   

\newcommand{\bR}{\mathbf{R}}

\newcommand{\Bu}{\boldsymbol{u}}

\newcommand{\HYP}{\mathbb{H}^3}
\newcommand{\HYN}{\mathbb{H}^n}

\begin{document}
\pagestyle{myheadings}
\markboth{\centerline{Mikl\'os Eper and Jen\H o Szirmai}}
{Coverings with horo- and hyperballs \dots}
\title
{Coverings with horo- and hyperballs generated by simply truncated orthoschemes}

\author{\normalsize{Mikl\'os Eper and Jen\H o Szirmai} \\
\normalsize Budapest University of Technology and \\
\normalsize Economics Institute of Mathematics, \\
\normalsize Department of Geometry \\
\date{\normalsize{\today}}}

\maketitle


\begin{abstract}

After having investigated the packings derived by horo- and hyperballs related to simple frustum Coxeter orthoscheme tilings we consider the corresponding covering problems 
(briefly hyp-hor coverings) in $n$-dimensional hyperbolic spaces $\HYN$ ($n=2,3$).

We construct in the $2-$ and $3-$dimensional hyperbolic spaces hyp-hor coverings that
are generated by simply truncated Coxeter orthocheme tilings
and we determine their thinnest covering configurations and their densities.

We prove that in the hyperbolic plane ($n=2$) the density of the above thinnest hyp-hor covering arbitrarily approximate
the universal lower bound of the hypercycle or horocycle covering density $\frac{\sqrt{12}}{\pi}$ and
in $\HYP$ the optimal configuration belongs to the $\{7,3,6\}$ Coxeter tiling with density $\approx 1.27297$ that is less than the previously known famous
horosphere covering density $1.280$ due to L.~Fejes T\'oth and K.~B\"or\"oczky.

Moreover, we study the hyp-hor coverings in
truncated orthosche\-mes $\{p,3,6\}$ $(6< p < 7, ~ p\in \mathbb{R})$ whose
density function attains its minimum at parameter $p\approx 6.45962$ with  density $\approx 1.26885$. That means that this
locally optimal hyp-hor configuration provide smaller covering density than the former determined $\approx 1.27297$ 
but this hyp-hor packing configuration can not be extended to the entirety of hyperbolic space $\mathbb{H}^3$.
\end{abstract}

\newtheorem{theorem}{Theorem}[section]
\newtheorem{corollary}[theorem]{Corollary}
\newtheorem{conjecture}{Conjecture}[section]
\newtheorem{lemma}[theorem]{Lemma}
\newtheorem{exmple}[theorem]{Example}
\newtheorem{defn}[theorem]{Definition}
\newtheorem{rmrk}[theorem]{Remark}
\newenvironment{definition}{\begin{defn}\normalfont}{\end{defn}}
\newenvironment{remark}{\begin{rmrk}\normalfont}{\end{rmrk}}
\newenvironment{example}{\begin{exmple}\normalfont}{\end{exmple}}
\newenvironment{acknowledgement}{Acknowledgement}


\section{Introduction}
The packing and covering problems with solely horo- or hyperballs (horo- or hypespheres) are intensively investigated
in earlier works in $n$-dimensional $(n\ge2)$ hyperbolic space $\HYN$.

In $n$-dimensional hyperbolic space $\mathbb{H}^n$ $(n\ge2)$ there are $3$ kinds
of ''balls (spheres)": the classical balls (spheres), horoballs (horospheres) and hyperballs (hyperspheres).

In this paper we consider the coverings with horo- and hyperballs and their densities in $2$- and $3$-dimensional hyperbolic space where the coverings are derived from  
simply truncated Coxeter orthoscheme tilings. 

A Coxeter simplex is an $n$-dimensional simplex in $X\in\{\mathbb{S}^n, \mathbb{H}^n, \mathbb{E}^n\}$ with dihedral angles either submultiples of $\pi$ or zero. 
The group generated by reflections on the sides of a Coxeter simplex is called a Coxeter simplex reflection group. 
Such reflections determine a discrete group of isometries of $X$ with the Coxeter simplex as its fundamental domain; 
hence such groups generate a tessellation of $X$. 

First we shortly survey the previous results related to this topic.
\begin{enumerate}
\item
{\bf On horoball packings and coverings}

In the case of periodic ball or horoball packings and coverings, the local density defined e.g. in \cite{B78} can be extended to the entire hyperbolic space. This local density
is related to the simplicial density function that we generalized in \cite{Sz12} and \cite{Sz12-2}.
In this paper we will use such definition of covering density.  

In the $n$-dimensional space $X \in\{ \mathbb{E}^n, \mathbb{S}^n, \mathbb{H}^n \}$ of constant curvature
$(n \geq 2)$, define the simplicial density function $d_n(r)$ to be the density of $n+1$ spheres
of radius $r$ mutually touching one another with respect to the regular simplex spanned by the centers of the spheres. L.~Fejes T\'oth and H.~S.~M.~Coxeter
conjectured that the packing density of balls of radius $r$ in $X$ cannot exceed $d_n(r)$.
Rogers \cite{Ro64} proved this conjecture in Euclidean space $\mathbb{E}^n$.
The $2$-dimensional spherical case was settled by L.~Fejes T\'oth \cite{FTL}, and B\"or\"oczky \cite{B78}, who proved the following extension:
\begin{theorem}[K.~B\"or\"oczky]
In an $n$-dimensional space of constant curvature, consider a packing of spheres of radius $r$.
In the case of spherical space, assume that $r<\frac{\pi}{4}$.
Then the density of each sphere in its Dirichlet--Voronoi cell cannot exceed the density of $n+1$ spheres of radius $r$ mutually
touching one another with respect to the simplex spanned by their centers.
\end{theorem}
In hyperbolic space $\HYP$, 
the monotonicity of $d_3(r)$ was proved by B\"or\"oczky and Florian
in \cite{B--F64}. 

This upper bound for packing density in hyperbolic space $\mathbb{H}^3$ is $\approx 0.85327$, 
which is not realized by packing regular balls. However, it is attained by a horoball packing of
$\overline{\mathbb{H}}^3$ where the ideal centers of horoballs lie on the
absolute figure of $\overline{\mathbb{H}}^3$;
for example, they may lie at the vertices of the ideal regular
simplex tiling with Coxeter-Schl\"afli symbol $\{3,3,6\}$. From this regular ideal tetrahedron tiling can be derived the known least dense ball or horoball covering configuration
(see \cite{FTL})with density $\approx 1.280$.

In \cite{KSz} we proved that the optimal ball packing arrangement in $\overline{\mathbb{H}}^3$ mentioned above is not unique. 
We gave several new examples of horoball packing arrangements based on totally asymptotic Coxeter tilings that yield the B\"or\"oczky--Florian upper bound \cite{B--F64}.
 
Furthermore, in \cite{Sz12}, \cite{Sz12-2} we found that 
by allowing horoballs of different types at each vertex of a totally asymptotic simplex and generalizing 
the simplicial density function to $\mathbb{H}^n$ for $(n \ge 2)$,
the B\"or\"oczky-type density 
upper bound is no longer valid for the fully asymptotic simplices for $n \ge 3$. 
For example, in $\overline{\mathbb{H}}^4$ the locally optimal packing density is $\approx 0.77038$, higher than the B\"or\"oczky-type density upper bound of $\approx 0.73046$. 
However these ball packing configurations are only locally optimal and cannot be extended to the entirety of the
hyperbolic spaces $\overline{\mathbb{H}}^n$. Further open problems and conjectures on $4$-dimensional hyperbolic packings are discussed in \cite{G--K--K}. 
Using horoball packings in $\mathbb{H}^4$, allowing horoballs of different types,
we found seven counterexamples (realized by allowing up to three horoball types) 
to one of L. Fejes T\'oth's conjectures stated in his foundational book Regular Figures.

In \cite{KSz19} and \cite{KSz19-1} we continued our investigations 
of ball packings, in hyperbolic spaces of dimensions $n=5\dots9$. 
Using horoball packings, allowing horoballs of different types when applicable,
we found several interesting and dense packing configuratons with respect to the Coxeter simplex cells.

The second-named author has several additional results on globally and locally optimal ball packings 
in the eight Thurston geomerties arising from Thurston's geometrization conjecture 
see e.g. \cite{Sz07-2}, \cite{Sz14-1}. 
 
\item
{\bf{On hyperball packings and coverings}}

In hyperbolic plane $\mathbb{H}^2$ the universal upper bound of the congruent hypercycle packing density is $\frac{3}{\pi}$,
proved by I.~Vermes in \cite{V79}. He initiated this topic and determined also the universal lower bound of the congruent hypercycle covering density, in \cite{V81},
equal to $\frac{\sqrt{12}}{\pi}$.

In \cite{Sz06-1} and \cite{Sz06-2} we have analysed the regular prism tilings (simple truncated Coxeter orthoscheme tilings) and the corresponding optimal hyperball packings in
$\mathbb{H}^n$ $(n=3,4,5)$. Recently (to the best of author's knowledge) these 
have been the densest packings with congruent hyperballs.

In \cite{Sz13-4} we studied the $n$-dimensional hyperbolic regular prism honeycombs
and the corresponding coverings by congruent hyperballs and we determined their least dense covering.
Furthermore, we formulated conjectures for the candidates of the least dense 
covering by congruent hyperballs in the 3- and 5-dimensional hyperbolic space. 

In \cite{Sz17-1} we discussed congruent and non-congruent hyperball packings to the truncated regular tetrahedron tilings.
These are derived from the truncated Coxeter simplex tilings 
$\{3,3,p\}$ $(7\le p \in \mathbb{N})$ and $\{3,3,3,3,5\}$
in $3$- and $5$-dimensional hyperbolic space, respectively.
We determined the densest packing arrangement and its density
with congruent hyperballs in $\mathbb{H}^5$ and determined the smallest density upper bounds of
non-congruent hyperball packings generated by the above tilings.

In \cite{Sz17} we deal with such packings by horo- and hyperballs (briefly hyp-hor packings) in $\HYN$
($n=2,3$).

In \cite{Sz14-2} we studied a large class of hyperball packings in $\HYP$
that can be derived from truncated tetrahedron tilings.
We proved that if the truncated tetrahedron is regular $\{3,3,p\}$, but we allow also 
$6< p \in \mathbb{R}$,
then the density
of the locally densest packing is $\approx 0.86338$. This is larger than the B\"or\"oczky-Florian density upper bound
but our locally optimal hyperball packing configuration cannot be extended to the entirety of
$\mathbb{H}^3$. However, we described a hyperball packing construction,
by the regular truncated 
tetrahedron tiling under the extended Coxeter group $\{3, 3, 7\}$ with maximal 
density $\approx 0.82251$.

{In \cite{Sz17-2} we developed a decomposition algorithm that for each saturated hyperball packing provides a decomposition of $\HYP$
into truncated tetrahedra. Therefore, in order to get a density upper bound for hyperball packings, it is sufficient to determine
the density upper bound of hyperball packings in truncated simplices.}

In \cite{Sz19} we proved, that the density upper bound of the saturated 
congruent hyperball packings, related to corresponding truncated tetrahedron cells is locally realized in a regular truncated 
tetrahedon with density $\approx  0.86338$.
Furthermore, we proved that the density of locally optimal congruent hyperball arrangement 
in regular truncated tetrahedron is not monotonically increasing function of the height 
of corresponding optimal hyperball, 
contrary to the ball and horoball packings.

In \cite{Sz19-4}, we considered hyperball packings related to truncated regular cube and octahedron tilings that are derived from the Coxeter truncated orthoscheme tilings
$\{4,3,p\}$ $(6< p \in \mathbb{N})$ and $\{3,4,p\}$ $(4 < p \in \mathbb{N})$
in hyperbolic space $\HYP$. If we allow $p \in \mathbb{R}$ as well, then the 
locally densest (non-congruent half) hyperball configuration belongs to the 
truncated cube with density
$\approx 0.86145$. This is larger than the B\"or\"oczky-Florian density upper bound for balls and horoballs.
But our locally optimal non-congruent hyperball packing configuration 
cannot be extended to the entire $\mathbb{H}^3$. We determined the extendable 
densest non-congruent hyperball packing arrangement related to the truncated cube tiling $\{4,3,7\}$
with density $\approx 0.84931$.

In \cite{Sz19-3} we studied congruent and non-congruent hyperball packings generated by doubly truncated Coxeter orthoscheme tilings in the $3$-di\-men\-sional hyperbolic space.
We proved that the densest congruent hyperball packing belongs to the Coxeter orthoscheme tiling of parameter $\{7,3,7\}$ with density $\approx 0.81335$. 
This density is equal -- by our conjecture -- with the upper bound density of the corresponding non-congruent hyperball arrangements.
\end{enumerate}
{\it In this paper} we deal with the coverings with horo- and hyperballs (briefly hyp-hor coverings) in the $n$-dimensional hyperbolic spaces $\HYN$
($n=2,3$) which form a new class of the classical covering problems.

We construct in the $2-$ and $3-$dimensional hyperbolic spaces hyp-hor coverings that
are generated by complete Coxeter tilings of degree $1$ i.e. the fundamental domains of these tilings are
simple frustum orthoschemes with a principal vertex lying on the absolute quadric 
and the other principal vertex is outer point.
We determine their thinnest covering configurations and their densities.
These considered Coxeter tilings exist in the $2-$, $3-$ and $5-$dimensional
hyperbolic spaces (see \cite{IH90}) and have given by their Coxeter-Schl\"afli graph in Fig.~1.
\begin{figure}[ht]
\centering
\includegraphics[width=12cm]{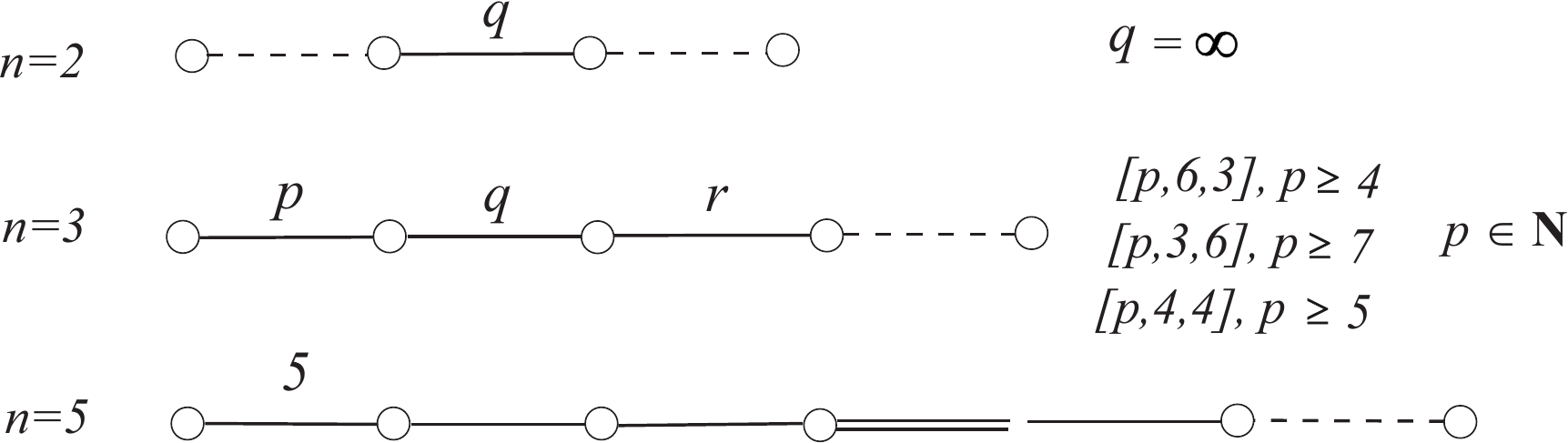}
\caption{Coxeter-Schl\"afli graph of Coxeter tilings of degree 1.}
\label{}
\end{figure}
We prove that in the hyperbolic plane $n=2$ the density of the above hyp-hor coverings arbitrarily approximate
the universal upper bound of the hypercycle or horocycle packing density $\frac{\sqrt{12}}{\pi}$ and
in $\HYP$ the thinnest hyp-hor configuration belongs to the $\{7,3,6\}$ Coxeter tiling with density $\approx 1.27297$.

Moreover, we consider the hyp-hor coverings in
truncated orthosche\-mes $\{p,3,6\}$ $(6< p < 7, ~ p\in \bR)$.
Its density function is attained its minimum for parameter $p \approx 6.45962$, and 
the corresponding minimal covering density is $\approx 1.26885$ less than $\approx 1.280$. That means that this
locally optimal hyp-hor configurations provide less densities that the previously known Fejes T\'oth-B\"or\"oczky-Florian covering density for ball and
horoball packings but this hyp-hor covering configurations can not be extended to the entirety of hyperbolic space $\mathbb{H}^3$.
\section{Basic notions}
For $\mathbb{H}^n$ we use the projective model in the Lorentz space $\mathbb{E}^{1,n}$ of signature $(1,n)$,
i.e.~$\mathbb{E}^{1,n}$ denotes the real vector space $\mathbf{V}^{n+1}$ equipped with the bilinear
form of signature $(1,n)$:~
$
\langle ~ \mathbf{x},~\mathbf{y} \rangle = -x^0y^0+x^1y^1+ \dots + x^n y^n
$
where the non-zero vectors
$
\mathbf{x}=(x^0,x^1,\dots,x^n)\in\mathbf{V}^{n+1} \ \  \text{and} \ \ \mathbf{y}=(y^0,y^1,\dots,y^n)\in\mathbf{V}^{n+1},
$
are determined up to real factors, for representing points of $\mathcal{P}^n(\mathbb{R})$. Then, $\mathbb{H}^n$ can be interpreted
as the interior of the quadric
$
Q=\{[\mathbf{x}]\in\mathcal{P}^n | \langle ~ \mathbf{x},~\mathbf{x} \rangle =0 \}=:\partial \mathbb{H}^n
$
in the real projective space $\mathcal{P}^n(\mathbf{V}^{n+1},
\mbox{\boldmath$V$}\!_{n+1})$.

The points of the boundary $\partial \mathbb{H}^n $ in $\mathcal{P}^n$
are called points at infinity of $\mathbb{H}^n $, the points lying outside $\partial \mathbb{H}^n $
are said to be outer points of $\mathbb{H}^n $ relative to $Q$. Let $P([\mathbf{x}]) \in \mathcal{P}^n$, a point
$[\mathbf{y}] \in \mathcal{P}^n$ is said to be conjugate to $[\mathbf{x}]$ relative to $Q$ if
$\langle ~ \mathbf{x},~\mathbf{y} \rangle =0$ holds. The set of all points which are conjugate to $P([\mathbf{x}])$
form a projective (polar) hyperplane 
$
{\boldsymbol{x}}=pol(\mathbf{x}):=\{[\mathbf{y}]\in\mathcal{P}^n | \langle ~ \mathbf{x},~\mathbf{y} \rangle =0 \}.
$
Thus the quadric $Q$ induces a bijection
(linear polarity $\mathbf{V}^{n+1} \rightarrow
\mbox{\boldmath$V$}\!_{n+1})$
from the points of $\mathcal{P}^n$
onto its hyperplanes.

The distance $s$ of two proper points
$[\mathbf{x}]$ and $[\mathbf{y}]$ is calculated by the formula:
\begin{equation}
\cosh{\frac{s}{k}}=\frac{-\langle ~ \mathbf{x},~\mathbf{y} \rangle }{\sqrt{\langle ~ \mathbf{x},~\mathbf{x} \rangle
\langle ~ \mathbf{y},~\mathbf{y} \rangle }} . \tag{2.1}
\end{equation}
\subsection{Complete orthoschemes}
A $n$-dimensional tiling $\mathcal{P}$ (or solid tessellation, honeycomb) is an infinite set of
congruent polyhedra (polytopes) that fit together to fill all space $(\mathbb{H}^n~ (n \geqq 2))$ exactly once,
so that every face of each polyhedron (polytope) belongs to another polyhedron as well.
At present the cells are congruent orthoschemes (see \cite{K89}).

Geometrically, complete orthoschemes of degree $d$ can be described as follows:
\begin{enumerate}
\item
For $d=0$, they coincide with the class of classical orthoschemes introduced by
{{Schl\"afli}}.
The initial and final vertices, $A_0$ and $A_n$ of the orthogonal edge-path
$A_iA_{i+1},~ i=0,\dots,n-1$, are called principal vertices of the orthoscheme.
\item
A complete orthoscheme of degree $d=1$ can be interpreted as an
orthoscheme with one outer principal vertex, say $A_n$, which is truncated by
its polar plane $pol(A_n)$ (see Fig.~2 and 3). In this case the orthoscheme is called simply truncated with
outer vertex $A_n$.
\item
A complete orthoscheme of degree $d=2$ can be interpreted as an
orthoscheme with two outer principal vertices, $A_0,~A_n$, which is truncated by
its polar hyperplanes $pol(A_0)$ and $pol(A_n)$. In this case the orthoscheme is called doubly
truncated. We distinguish two different types of orthoschemes but I
will not enter into the details (see \cite{K89}).
\end{enumerate}

In general the complete Coxeter orthoschemes were classified by {{Im Hof}} in
\cite{IH90} by generalizing the method of {{Coxeter}} and {{B\"ohm}}, who
showed that they exist only for dimensions $\leq 9$. From this classification it follows, that the complete
orthoschemes of degree $d=1$ exist up to 5 dimensions. 

In this paper we consider the orthoschemes of degree 1 where the initial vertex $A_0$ lies on the
absolute quadric $Q$. These orthoschemes and the corresponding Coxeter tilings exist in the $2$-, $3-$ and $5-$dimensional hyperbolic spaces and
are characterized by their Coxeter-Schl\"afli symbols and graphs (see Fig.~1).

In $n$-dimensional hyperbolic space $\mathbb{H}^n$ $(n \ge 2)$
it can be seen that if $\mathcal{S}$ is a complete
orthoscheme of degree $d=1$ (with vertices $A_0A_1A_2 \dots A_{n-1}$ $P_0P_1P_2 \dots P_{n-1}$) a simply frustum orthoscheme (here $A_n$ is a outer vertex of
$\mathbb{H}^n$ then the points $P_0,P_1,P_2,\dots,P_{n-1}$ lie on the polar hyperplane $\pi$ of $A_n$).

We consider the images of $\mathcal{S}$ under reflections on its side facets.
The union of these $n$-dimensional orthoschames (having the common $\pi$ hyperplane) forms an infinite polyhedron denoted by $\mathcal{G}$.
$\mathcal{G}$ and its images under reflections on its ,,cover facets" fill hyperbolic
space $\mathbb{H}^n$ without overlap and generate $n$-dimensional tilings $\mathcal{T}$.

{\it The constant $k =\sqrt{\frac{-1}{K}}$ is the natural length unit in
$\mathbb{H}^n$. $K$ will be the constant negative sectional curvature. In the following we assume that $k=1$.}
\subsection{Volumes of the $n$-dimensional \\ Coxeter orthoschemes}

\begin{enumerate}
\item $2$-dimensional hyperbolic space $\mathbb{H}^2$

In the hyperbolic plane a simple frustum orthoscheme is a Lambert quadrilateral with exactly three right angles and its fourth angle is acute
$\frac{\pi}{q}$ ($q \ge 3$) (see Fig.~1 and 3). In our case the Lambert quadrilateral has a vertex at the infinity i.e. the angle at this vertex is $0$.
Its area can be determined by the well-known defect formula of hyperbolic triangles:
\begin{equation}
Vol_2(\mathcal{S})=\frac{\pi}{2}. \tag{2.2}
\end{equation}

\item $3$-dimensional hyperbolic space $\HYP$:

{Our polyhedron $A_0A_1A_2P_0P_1P_2$ is a simple frustum orthoscheme with
outer vertex $A_3$ (see Fig.~5.a) whose volume can be calculated by the following theorem of R.~Kellerhals
\cite{K89}:}
\begin{theorem} The volume of a three-dimensional hyperbolic
complete ortho\-scheme (except Lambert cube cases) $\mathcal{S}$
is expressed with the essential angles $\alpha_{01},\alpha_{12},\alpha_{23}, \ (0 \le \alpha_{ij} \le \frac{\pi}{2})$
(Fig.~1 and 2) in the following form:

\begin{align}
&Vol_3(\mathcal{S})=\frac{1}{4} \{ \mathcal{L}(\alpha_{01}+\theta)-
\mathcal{L}(\alpha_{01}-\theta)+\mathcal{L}(\frac{\pi}{2}+\alpha_{12}-\theta)+ \notag \\
&+\mathcal{L}(\frac{\pi}{2}-\alpha_{12}-\theta)+\mathcal{L}(\alpha_{23}+\theta)-
\mathcal{L}(\alpha_{23}-\theta)+2\mathcal{L}(\frac{\pi}{2}-\theta) \}, \tag{2.3}
\end{align}
where $\theta \in [0,\frac{\pi}{2})$ is defined by the following formula:
$$
\tan(\theta)=\frac{\sqrt{ \cos^2{\alpha_{12}}-\sin^2{\alpha_{01}} \sin^2{\alpha_{23}
}}} {\cos{\alpha_{01}}\cos{\alpha_{23}}}
$$
and where $\mathcal{L}(x):=-\int\limits_0^x \log \vert {2\sin{t}} \vert dt$ \ denotes the
Lobachevsky function.
\end{theorem}
For our prism tilings $\mathcal{T}_{pqr}$ we have:~
$\alpha_{01}=\frac{\pi}{p}, \ \ \alpha_{12}=\frac{\pi}{q}, \ \
\alpha_{23}=\frac{\pi}{r}$ .
\end{enumerate}
\subsection{On hyperballs}
The equidistant surface (or hypersphere) is a quadratic surface that lies at a constant distance
from a plane in both halfspaces. The infinite body of the hypersphere is called a hyperball.
The $n$-dimensional {\it half-hypersphere } $(n=2,3)$ with distance $h$ to a hyperplane $\pi$
is denoted by $\mathcal{H}_n^h$.
The volume of a bounded hyperball piece $\mathcal{H}_n^h(\mathcal{A}_{n-1})$
bounded by an $(n-1)$-polytope $\mathcal{A}_{n-1} \subset \pi$, $\mathcal{H}_n^h$ and by
hyperplanes orthogonal to $\pi$ derived from the facets of $\mathcal{A}_{n-1}$ can be determined by the
formulas (2.4) and (2.5) that follow from the suitable extension of the classical method of {{J.~Bolyai}} (\cite{B91}):
\begin{equation}
Vol_2(\mathcal{H}_2^h(\mathcal{A}_1))=Vol_1(\mathcal{A}_{1}) \sinh{(h)}, \tag{2.4}
\end{equation}
\begin{equation}
Vol_3(\mathcal{H}_3^h(\mathcal{A}_2))=\frac{1}{4}Vol_2(\mathcal{A}_{2})\left[\sinh{(2h)}+
2 h \right], \tag{2.5}
\end{equation}
where the volume of the hyperbolic $(n-1)$-polytope $\mathcal{A}_{n-1}$ lying in the plane
$\pi$ is $Vol_{n-1}(\mathcal{A}_{n-1})$.
\subsection{On horoballs}
A horosphere in $\mathbb{H}^n$ $(n \ge 2)$ is a
hyperbolic $n$-sphere with infinite radius centered
at an ideal point on $\partial \mathbb{H}^n$. Equivalently, a horosphere is an $(n-1)$-surface orthogonal to
the set of parallel straight lines passing through a point of the absolute quadratic surface.
A horoball is a horosphere together with its interior.

We consider the usual Beltrami-Cayley-Klein ball model of $\mathbb{H}^n$
centered at $O(1,0,0,$ $\dots, 0)$.
The equation of a horosphere with center
$T_0(1,0,\dots,1)$ passing through the point $S(1,0,\dots,s)$ is derived from the equation of the
the absolute sphere $-x^0 x^0 +x^1 x^1+x^2 x^2+\dots + x^n x^n = 0$, and the plane $x^0-x^n=0$ tangent to the absolute sphere at $T_0$.
The general equation of the horosphere in cartesian coordinates is the following:
\begin{equation}
\label{eqn:horosphere1}
\frac{2 \left(\sum_{i=1}^{n-1} h_i^2 \right)}{1-s}+\frac{4 \left(h_n-\frac{s+1}{2}\right)^2}{(1-s)^2}=1. \tag{2.6}
\end{equation}

{\it In $n$-dimensional hyperbolic space any two horoballs are congruent in the classical sense.
However, it is often useful to distinguish between certain horoballs of a packing.
We use the notion of horoball type with respect to the packing as introduced in \cite{Sz12-2}.}

The intrinsic geometry of a horosphere is Euclidean,
so the $(n-1)$-dimensional volume $\mathcal{A}$ of a polyhedron $A$ on the
surface of the horosphere can be calculated as in $\mathbb{E}^{n-1}$.
The volume of the horoball piece $\mathcal{H}(A)$ determined by $A$ and
the aggregate of axes
drawn from $A$ to the center of the horoball is (\cite{B91})
\begin{equation}
\label{eq:bolyai}
Vol(\mathcal{H}(A)) = \frac{1}{n-1}\mathcal{A}. \tag{2.7}
\end{equation}
\section{Hyp-hor coverings in hyperbolic plane}
We consider the usual Beltrami-Cayley-Klein ball modell of $\mathbb{H}^2$ centered at $O(1,0,0)$ with a given vector basis ${\bf{e}}_i$ $(i=0,1,2)$ and set the 2-dimensional Coxeter orthoscheme $A_0A_1A_2$ in this coordinate system with coordinates $A_0(1,0,1); A_1(1,0,0); A_2(1,\frac{1}{a},0)$. Here the initial principal vertex of the orthoscheme $A_0$ is lying on the absolute quadric $Q$ and the other principal vertex $A_2$ is an outer point of the model, so $ 0 < a < 1, a \in \mathbb{R}$.

The polar line of the outer vertex $A_2$ is $\pi={\Bu}_2(1,-\frac{1}{a},0)^T$. By the truncation of the orthoscheme $A_0A_1A_2$ by the polar line 
$\pi$ we get the Lambert quadrilateral $A_0A_1P_1P_0$ (see Fig. 2), where the further vertices are: $\pi\cap A_0A_2 = P_0(1,a,1-a^2); \pi \cap A_1A_2=P_1(1,a,0)$. Its images under reflections on its sides fill hyperbolic plane $\mathbb{H}^2$ without overlap, hence we get the previously described 2-dimensional Coxeter tilings, given by the Coxeter symbol $\left[ \infty \right]$ (see Fig. 1). The tilings contain the free parameter $a$, so we denote the tilings by $\mathcal{T}_a$, and the Lambert quadrilaterals $A_0A_1P_0P_1$ by $\mathcal{F}_a$, which serve as the fundamental domain of the above tilings.

\begin{figure}[ht]
\centering
\includegraphics[width=6.5cm]{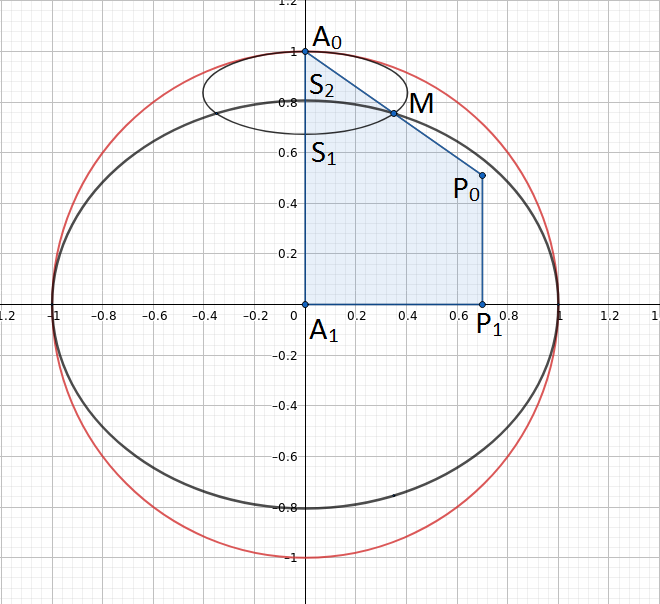} 
\includegraphics[width=6cm]{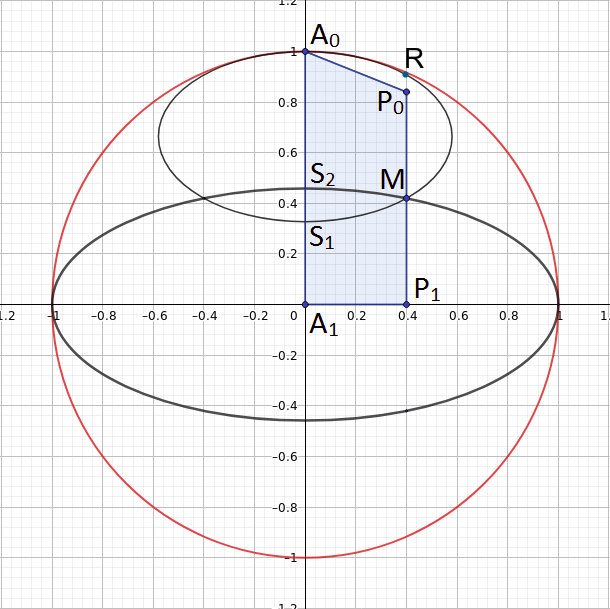}

a) \hspace{5cm} b)
\caption{a)~$\mathcal{C}_a^1$-type hyp-hor covering at present $a=0.7,t=0.5$ ~b)~$\mathcal{C}_a^2$-type hyp-hor covering at present $a=0.4,t=0.5$}
\end{figure}

We construct hyp-hor coverings to $\mathcal{F}_a$, by the follows:
\begin{itemize}
\item[1. ]The center of the horocycle can only be the vertex $A_0$. Let the intersection of the horocycle with $A_0A_1$ line $S_1(1,0,s_1)$ $(-1 < s_1 < 1)$ 
and with $A_0P_0$ line $T(1,ta,1-ta^2)$ $(0 < t < \frac{2}{1+a^2})$. We denote by $\mathfrak{H}_a(t)$ the horocycle-piece determined by points $A_0,S_1,T$ (see Fig. 2).
\item[2. ] Let $A_1P_1$ be the base straight line of a hypercycle and $M$ the intersection point of the horo- and hypercycle lies on the $A_0P_0$ or $P_0P_1$ side 
of $\mathcal{F}_a$ (see Fig. 2).
\item[3. ] Let the intersection of the hypercycle with the positive segment of $A_0A_1$ line $S_2(1,0,s_2)$ $(0 < s_2 < 1)$ and with $P_0P_1$ line $R(1,a,r)$ $(0 < r < \sqrt{1-a^2})$. We denote by $\mathcal{H}_a(t)$ the hypercycle-piece settled by points $P_1,R,S_2,A_1$ (see Fig. 2).
\end{itemize}
We can see, that if the horo- and hypercycles satisfy the above requirements, than they cover $\mathcal{F}_a$. 
Thus the images of $\mathfrak{H}_a(t)$ and $\mathcal{H}_a(t)$ under reflection on the sides of $\mathcal{F}_a$ provide a hyp-hor covering of hyperbolic plane $\mathbb{H}^2$. 
The fundamental domain $\mathcal{F}_a$ (i.e. parameter $a$) and point $M$ (i.e. parameter $t$) determine the covering. We distinguish two main types of hyp-hor coverings, denoted by $\mathcal{C}^1_a(t)$ if $M \in A_0P_0$ and by $\mathcal{C}^2_a(t)$ if $M \in P_0P_1$ (see Fig. 2).

\begin{defn}
The density of the above hyp-hor coverings $\mathcal{C}^i_a(t)$ $(i=1,2)$ are:
\begin{equation*}
\delta(\mathcal{C}^i_a(t))=\frac{Vol(\mathcal{H}_a(t))+Vol(\mathfrak{H}_a(t))}{Vol(\mathcal{F}_a)}
\end{equation*}
\end{defn}

It is obvious, that if the point $M$ lies on the perimeter of $\mathcal{F}_a$, the density of the covering is smaller, than it lies out of $\mathcal{F}_a$. Thus we get the coverings with minimal densities in the above two main cases. 
\subsection{The densities of coverings $\mathcal{C}^1_a(t).$}
In this case $M \in A_0P_0$ is the intersection point of the cycles, so $M=(1,ta,1-ta^2)$ $(0< t \leq 1)$. 
The coordinates of $S_1$ can be expressed using (2.6) and the distance of $M$ and $S_1$ can be calculated by (2.1), thus we can determine the volume of $\mathfrak{H}_a(t)$ 
by formula (2.6). The length of $A_1P_1$ and the distance of $M$ and the $x-$axis can be calculated also by (2.1), thus we can determine the volume of $\mathcal{H}_a(t)$ by formula (2.4). We obtain by Definition 3.1, that the density of $\mathcal{C}^1_a(t)$ can be expressed by the following formula:

\begin{equation*}
\delta(\mathcal{C}_a^1(t))=\frac{\textmd{arccosh}\left(\frac{1}{\sqrt{1-a^2}}\right)\frac{1-ta^2}{a\sqrt{2t-t^2-a^2t^2}}+2\sinh\left(\frac{1}{2}\textmd{arccosh}\left(\frac{2ta^2+t-4}{2t-4+2ta^2}\right)\right)}{\frac{\pi}{2}}
\end{equation*}

where $0 < a < 1$, $0 < t \leq 1$.

\begin{theorem}
Analysing the above density formula we obtain that
\begin{equation*}
\lim_{a\rightarrow 0}\left(\mathcal{C}_a^1\left(\frac{1}{2}\right)\right)=\frac{\sqrt{12}}{\pi}
\end{equation*}
and $\left(\mathcal{C}_a^1\left(\frac{1}{2}\right)\right)<\frac{\sqrt{12}}{\pi}$ for parameter $0 < a < 1$ (see Fig. 3a).
That means, that in hyperbolic plane $\mathbb{H}^2$ the universal lower bound density of ball and horoball coverings can be arbitrary accurate approximate with the densities $\delta\left(\left(\mathcal{C}_a^1\left(\frac{1}{2}\right)\right)\right)$ of hyp-hor packings of type 1.
\end{theorem}
\begin{figure}[ht]
\centering
\includegraphics[width=6.5cm]{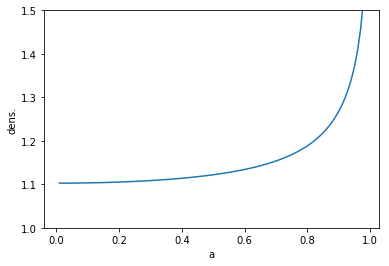} 
\includegraphics[width=6.5cm]{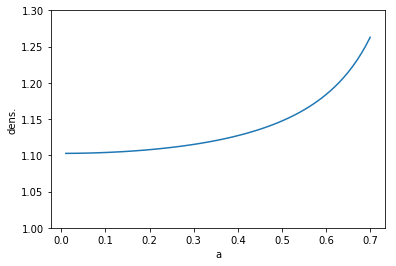}

a) \hspace{5cm} b)
\caption{a)~The density function of hyp-hor covering $\mathcal{C}_a^1$ in case $t=0.5$~b)~The density function of hyp-hor covering $\mathcal{C}_a^2$ in case $t\approx1.142$}
\end{figure}
\subsection{The densities of coverings $\mathcal{C}^2_a(t).$}
In this case $M \in P_0P_1$ the intersection point of the cycles, so the intersection point of the horocycle and line $A_0P_0$ is $(1,ta,1-ta^2)$ 
$(0< t < \frac{2}{1+2a^2-a^4})$, by the condition, that $M$ lies on the positive segment of $P_0P_1$. We get the volume of $\mathfrak{H}_a(t)$ just like in 
the previous section. The coordinates of $M$ and the $h_2$  length of $MP_1$ can be calculated by (2.6) and (2.1). 
We can determine the volume of $\mathcal{H}_a(t)$ by formula (2.4). We obtain by Definition 3.1, that the density of $\mathcal{C}^2_a(t)$ can be expressed by the following formula:

\begin{equation*}
\delta(\mathcal{C}_a^2(t))=\frac{\textmd{arccosh}\left(\frac{1}{\sqrt{1-a^2}}\right)\sinh{h_2}+2\sinh\left(\frac{1}{2}\textmd{arccosh}\left(\frac{2ta^2+t-4}{2t-4+2ta^2}\right)\right)}{\frac{\pi}{2}}
\end{equation*}

where $0 < a < 1$, $0 < t < \frac{2}{1+2a^2-a^4}$.

\begin{theorem}
Analysing the above density formula (using also numerical approximation methods) we obtain that it provides its minimum in case $t \approx 1.142$, $a\rightarrow 0$ (see Fig. 3b), and the minimum value is $\frac{\sqrt{12}}{\pi}$. That means, that in hyperbolic plane $\mathbb{H}^2$ the universal lower bound density of ball coverings can be arbitrary accurate approximate with the densities $\delta\left(\mathcal{C}_a^2\right)$ of hyp-hor packings of type 2.
\end{theorem}
\section{Hyp-hor coverings in hyperbolic space $\mathbb{H}^3$}
In the $3$-dimensional hyperbolic space $\HYP$ there are $3$ infinite series of the simple frustum Coxeter orthoschemes with vertex at the infinity, that are 
listed in Fig.~1, and characterized in Section 2.1. The Coxeter-Schl\"afli symbol of these orthoschemes are $\{p,q,r\}$, where $(q,r)=(3,6),(4,4),(6,3)$, and $p$ is an
appropriate integer parameter: $p\geq7$ if $(q,r)=(3,6)$, $p\geq5$ if $(q,r)=(4,4)$, $p\geq4$ if $(q,r)=(6,3)$. 
These conditions came from the geometry of the orthoschemes and can be computed by the inverse Coxeter-Schl\"afli matrix. We denote the orthoscheme by 
$\mathcal{F}^{(q,r)}_p$, and its vertices are denoted by $A_0$, $A_1$, $A_2$, $P_0$, $P_1$, $P_2$ (see Fig.~5.a). 

We consider the usual Beltrami-Cayley-Klein ball modell of $\mathbb{H}^3$ centred at $O(1,0,0,0)$ with a given vector basis ${\bf{e}}_i$ $(i=0,1,2,3)$ (see Section 2.1) 
and with the 3-dimensional complete Coxeter orthoscheme $A_0A_1A_2A_3$ which initial principal vertex $A_0$ is lying on the absolute quadric 
$Q$ and the other principal vertex $A_3$ is an outer point of the model. By the truncation of the orthoscheme with $\pi$ (the polar plane of $A_3$) we get 
the proper vertices $P_k[\mathbf{p}_k]=\pi\cap A_kA_3,(i=0,1,2)$, therefore ${\bf{p}}_k\thicksim c\cdot {\bf{a}}_3+{\bf{a}}_k$ for some $c\in\mathbb{R}$. $P_k[{\bf{p}}_k]$ lies on ${\boldsymbol{a}}^3=pol({\bf{a}}_3)$ if and only if ${\bf{p}}_k {\boldsymbol{a}}^3=0$, thus:
\begin{gather}
c\cdot {\bf{a}}_3{\boldsymbol{a}}^3+{\bf{a}}_k{\boldsymbol{a}}^3=0 \Leftrightarrow c=-\frac{{\bf{a}}_k{\boldsymbol{a}}^3}{{\bf{a}}_3{\boldsymbol{a}}^3} \tag{4.1}
\\ \Leftrightarrow {\bf{p}}_k \thicksim -\frac{{\bf{a}}_k{\boldsymbol{a}}^3}{{\bf{a}}_3{\boldsymbol{a}}^3}\cdot {\bf{a}}_3+{\bf{a}}_k \thicksim {\bf{a}}_k({\bf{a}}_3{\boldsymbol{a}}^3)+{\bf{a}}_3({\bf{a}}_k{\boldsymbol{a}}^3) \tag{4.2}
\end{gather}

We consider the Coxeter-Schl\"afli matrix $(c^{ij})$ of the orthoscheme, and its inverse $(h_{ij})$, where the elements of the matrices: 
$c^{ij}={\boldsymbol{a}}^i{\boldsymbol{a}}^j$, $h_{ij}={\bf{a}}_i{\bf{a}}_j$. The polar hyperplane of $A_3$ is 
${\boldsymbol{a}}^3$, thus $h_{k3}={\bf{a}}_k{\boldsymbol{a}}^3$, hence by (4.2) ${\bf{p}}_k={\bf{a}}_k h_{33}-{\bf{a}}_3 h_{k3}$.

We set the above simple frustum orthoscheme $\mathcal{F}^{(q,r)}_p$ in the usual coordinate system with vertices: $P_0(1,0,0,0)$, 
$P_1(1,0,y,0)$, $P_2(1,x,y,0)$, $A_0(1,0,0,1)$, $A_1(1,0,y,z_1)$, $A_2(1,x,y,z_2)$ (see Fig.~5.a). 
We get the following equations, using the formulas (2.1), (4.2) and $h_{ij}={\bf{a}}_i{\bf{a}}_j$:

\begin{gather*}
\cosh(d(P_0P_1))=\frac{h_{03}h_{13}-h_{01}h_{33}}{\sqrt{(h_{11}h_{33}-h_{13}^2)(h_{00}h_{33}-h_{03}^2)}}=\frac{1}{\sqrt{(-1)(-1+y^2)}}, \tag{4.3}
\\ \cosh(d(P_0P_2))=\frac{h_{03}h_{23}-h_{02}h_{33}}{\sqrt{(h_{22}h_{33}-h_{23}^2)(h_{00}h_{33}-h_{03}^2)}}=\frac{1}{\sqrt{(-1)(-1+y^2+x^2)}}, \tag{4.4}
\\ \cosh(d(A_1P_1))=\sqrt{\frac{h_{11}h_{33}-h_{13}^2}{h_{11}h_{33}}}=\frac{1-y^2}{\sqrt{(-1+y^2+z_1^2)(-1+y^2)}}, \tag{4.5}
\\ \cosh(d(A_2P_2))=\sqrt{\frac{h_{22}h_{33}-h_{23}^2}{h_{22}h_{33}}}=\frac{1-y^2-x^2}{\sqrt{(-1+y^2+x^2+z_2^2)(-1+y^2+x^2)}}. \tag{4.6}
\end{gather*} 
We can determine the coordinates $x,y,z_k,(k=1,2)$ by solving these equations, and the volume of the orthoschemes $\mathcal{F}^{(q,r)}_p$ by Theorem 2.1. 

The images of the above orthoscheme $\mathcal{F}^{(q,r)}_p$ under reflections on its facets fill the hyperbolic space $\mathbb{H}^3$ without overlap, 
so we get the Coxeter tiling $\mathcal{T}^{(q,r)}_p$ of $\mathbb{H}^3$ with fundamental domain $\mathcal{F}^{(q,r)}_p$.

We construct hyp-hor coverings to $\mathcal{F}^{(q,r)}_p$ using the following requirements:
\begin{itemize}
\item[1. ]The center of the horoball can only be the ideal vertex $A_0$. Let $S_1,T_1,Q_1$ the intersection points of the horoball with $A_0P_0,A_0A_2,A_0A_1$ 
lines. We denote by $\mathfrak{H}_p^{(q,r)}$ the horoball-piece determined by points $A_0,S_1,T_1,Q_1$ (see Fig.~5.b).
\item[2. ]$P_0P_1P_2$ plane can be the base hyperplane of a hyperball. 
Let $S_2,V_2,R_2$ be the intersection points of the hyperball with the line segments of $A_0P_0,A_1P_1,A_2P_2$. 
We denote by $\mathcal{H}_p^{(q,r)}$ the hyperball-piece bounded by the base hyperplane, the surface of the hyperball and the hyperplanes perpendicular to the base hyperplane 
derived from edges $P_0P_1,P_1P_2,P_2P_0$ (see Fig.~5.b).
\item[3. ]The intersection curve (which is a circle parallel with $[xy]$ plane in Euclidean sense) of the horo- and hyperball passes through one of the edges 
of the orthoscheme $A_0A_1,A_0A_2,A_1A_2,A_0P_0,A_1P_1,A_2P_2$ (see Fig.~5.a). 
\end{itemize}

We can see, that if the horo- and  hyperballs satisfy the above requirements, than they 
cover $\mathcal{F}^{(q,r)}_p$ if and only if they cover all the edges of $\mathcal{F}^{(q,r)}_p$. Hence, if a covering arrangement covers the edges 
of the orthoscheme, than the images of $\mathfrak{H}_p^{(q,r)}$ and $\mathcal{H}_p^{(q,r)}$ under reflection on the facets of $\mathcal{F}^{(q,r)}_p$ 
provide a hyp-hor covering of hyperbolic space $\mathbb{H}^3$, denoted by $\mathcal{C}_p^{(q,r)}$.

\begin{defn}
The density of the above hyp-hor coverings $\mathcal{C}^{(q,r)}_p$ is:
\begin{equation*}
\delta(\mathcal{C}^{(q,r)}_p)=\frac{Vol(\mathcal{H}_p^{(q,r)})+Vol(\mathfrak{H}_p^{(q,r)})}{Vol(\mathcal{F}^{(q,r)}_p)} \tag{4.7}
\end{equation*}
\end{defn}

It is obvious, that if the intersection curve passes through one of  the edge of $\mathcal{F}^{(q,r)}_p$, the density of the covering is smaller, 
than it goes out of $\mathcal{F}^{(q,r)}_p$. Thus we get the coverings with minimal densities if the above requiements hold. 
Based on the above, we have to distinguish and study six cases. 
\subsection{Non-covering cases}
\begin{itemize}
\item If the intersection curve of the balls passes through $A_0P_0$ (see Fig.~5.a), then the balls touch each other, thus the hyp-hor covering is obviously not realized.
\item If the intersection curve of the balls intersects the edge $A_0A_2$ (see Fig.~5.a), then we can parametrize their common point: 
$T(t)=(1,tx,ty,tz_2+(1-t)),t\in [0,1]$. By substituting this in the equation of the balls, we get the coordinates of $S_1\in P_0A_0, S_2 \in P_0A_0$ points. 
If the horoball covers $A_1$, we can determine the intersection points $U_1,U_2$ by 
solving the corresponding equations. By inspecting the $z$-coordinates of $U_i$ $(i=1,2)$ in the model, we can see, that $U_1$ is always higher than $U_2$,
which means (using the convexity of the ellipsoids) that they together do not cover the edge $A_1A_2$. If the hyperball covers $A_1$, we can determine 
the intersection points $Q_1,Q_2$ by solving the corresponding equations. By inspecting the $z$-coordinates of $Q_i$ $(i=1,2)$ in the model, we can see, 
that $Q_1$ is always higher than $Q_2$, which means (using the convexity of the ellipsoids) that they together do not cover the edge $A_0A_1$. Thus in 
this case the hyp-hor covering is not realized.
\item If the intersection curve of the balls contains a point of $A_1P_1$ edge (see Fig.~5.~a) then we can parametrize the intersection point $V$:
$V(v)=(1,0,y,vz_1),v\in [0,1]$. Very similarly to the above case, we can see, that if the horoball covers $A_2$, than the balls do not cover edge $A_2P_2$, 
and if the hyperball covers $A_2$, than the balls do not cover edge $A_1A_2$. Thus, in this case the hyp-hor covering is not realized. 
\end{itemize}
\subsection{Thinnest covering, if the intersection point lies on $A_0A_1$ edge}
In this case, $A_0A_1$ edge has a common point with the intersection curve of the balls (see Fig.~5.a), so we can parametrize 
the intersection point $Q$: $Q(q)=(1,0,qy,qz_1+(1-q)),q\in [0,1]$. By substituting this in the equation of the balls, we get the coordinates 
of $S_1,S_2 \in P_0A_0$. After that, we can determine the intersection points $T_1,T_2\in A_0A_2$ by solving the corresponding equations. 
We prove, that the balls cover the edges of the orthoscheme, so the hyp-hor covering is realized in this case. $P_0A_0A_1P_1$ is a $2$-dimensional Coxeter orthoscheme, 
thus $A_0A_1$ is covered as we have seen in Section 3. The hyperball covers $A_1$, and we can see, that the hyperbolic length of $A_1P_1$ edge is always bigger than the length 
of $A_2P_2$ edge, so the hyperball covers $A_2$, and because of its convexity $A_1P_1,A_2P_2,A_1A_2$ edges as well. By 
inspecting the $z$-coordinates of $S_i$ and $T_i$ $(i=1,2)$ in the model, we can see, that $S_2$ is always ``higher" than $S_1$ and $T_2$ is always ``higher" than $T_1$, 
which means (using the convexity of the ellipsoids) that they together cover the edges $A_0P_0$ and $A_0A_2$.

We know the coordinates of points $Q,T_i,S_i$ $(i=1,2)$, so we can determine the $Vol(\mathcal{H}_p^{(q,r)})$,  $Vol(\mathfrak{H}_p^{(q,r)})$ using (2.5), (2.7) 
and the density of the covering using (4.7), which depends on free parameter $q$. Analysing this density function we can compute the optimal densities (see Fig.~4.a). 
The results for tiling $\mathcal{T}_p^{(6,3)}$ (which provides the lowest density in this case) are summarized in the table below.
\\ 
\begin{center}
\begin{tabular}{|c|c|c|}
\hline
\multicolumn{1}{|c}{Type of tiling} & \multicolumn{1}{|c|}{$\delta_{min}$} & \multicolumn{1}{|c|}{$q$}\\
\hline
$\mathcal{T}_4^{(6,3)}$      &$1.3482413$ &$0.7369142$ \\
\hline
$\mathcal{T}_5^{(6,3)}$      &$1.4432379$ &$0.7655641$ \\
\hline
$\mathcal{T}_6^{(6,3)}$      &$1.5178400$ &$0.7814085$ \\
\hline
\end{tabular}
\end{center}
\begin{figure}[ht]
\centering
\includegraphics[width=6.5cm]{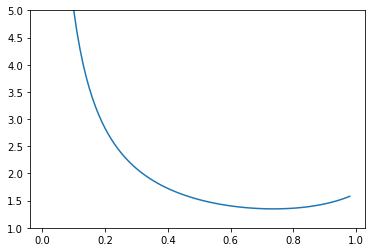} 
\includegraphics[width=6.5cm]{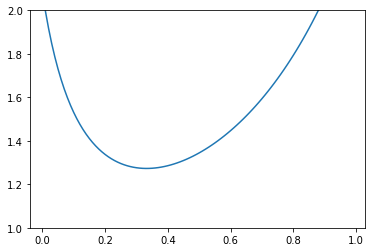}

a) \hspace{5cm} b)
\caption{a)~The density function $\delta(\mathcal{C}^{(6,3)}_4(q))$ ~b)~The density function $\delta(\mathcal{C}^{(3,6)}_7(u))$}
\end{figure}
\subsection{Thinnest covering, if the intersection point lies on $A_1A_2$ edge}
Now, the intersection curve of the balls passes through $A_1A_2$ (see Fig.~5.a), so we can parametrize the intersection point $U \in A_1A_2$: 
$U(u)=(1,ux,y,uz_2+(1-u)z_1),u\in [0,1]$. By substituting this in the equation of the balls, we get the coordinates of $S_1,S_2$ points. 
After that we can determine the intersection points $V_1,V_2$, $Q_1,Q_2$ and $T_1,T_2$ by solving the corresponding equations. 
We can prove similarly to the above case, that the balls cover the edges of the 
orthoscheme, so the hyp-hor covering is realized in this case. The horoball covers $A_0A_1$, the hyperball covers $A_2P_2$, and together they cover $A_1A_2$ (see Fig.~5.b). 
By inspecting the $z$-coordinates of $S_i$, $V_i$ and $T_i$ $(i=1,2)$ in the model, we can see in this case too, that the balls cover $A_0P_0,A_1P_1,A_0A_2$ edges (see Fig.~5.b).

We know points $Q_i,T_i,S_i$ $(i=1,2)$, so we can determine the $Vol(\mathcal{H}_p^{(q,r)})$,  $Vol(\mathfrak{H}_p^{(q,r)})$ using (2.5), (2.7) and the density of the 
covering using (4.7), which depends on free parameter $u$. Analysing this density function we can compute the optimal densities (see Fig. 4 b). 
The results for tiling $\mathcal{T}_p^{(3,6)}$ (which provides the lowest density in this case) are summarized in the next table.
\begin{figure}[ht]
\centering
\includegraphics[width=5cm]{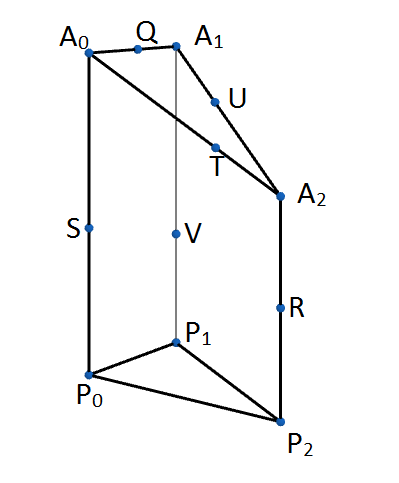} 
\includegraphics[width=8.5cm]{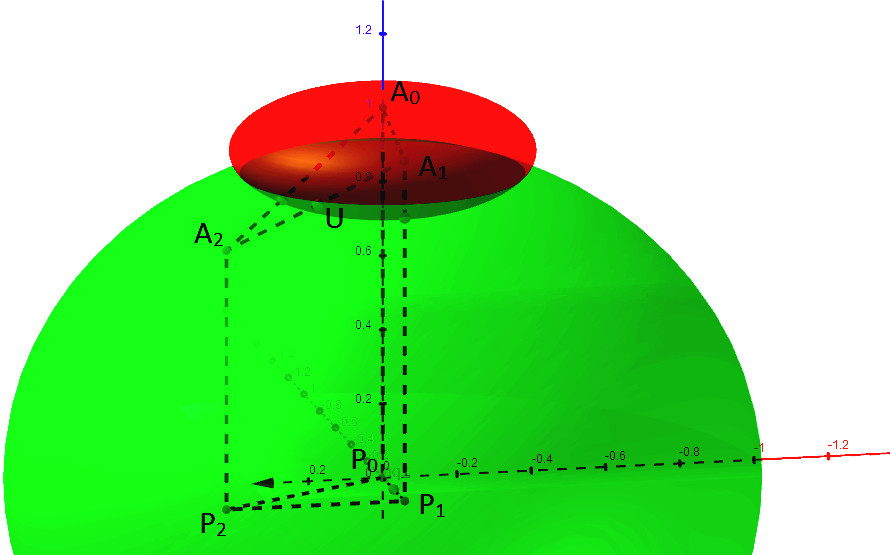}

a) \hspace{7cm} b)
\caption{a)~Simple truncated orthoscheme, and the intersection points of the balls, in the 6 cases ~b)~Hyp-hor covering of $\mathcal{F}_7^{(3,6)}$ with smallest density $\approx1.27297$}
\end{figure}
\\ 
\begin{center}
\begin{tabular}{|c|c|c|}
\hline
\multicolumn{1}{|c}{Type of tiling} & \multicolumn{1}{|c|}{$\delta_{min}$} & \multicolumn{1}{|c|}{$u$}\\
\hline
$\mathcal{T}_7^{(3,6)}$      &$1.27297329$ &$0.3324288$ \\ 
\hline
$\mathcal{T}_8^{(3,6)}$      &$1.288832$ &$0.3337034$ \\
\hline
$\mathcal{T}_9^{(3,6)}$      &$1.3065421$ &$0.3358650$ \\
\hline
\end{tabular}
\end{center}
\begin{rmrk} 
To any parameter $p$ $(6<p<7,p\in\mathbb{R})$ belongs a simple frustum orthoscheme $\mathcal{F}_p^{(3,6)}$ as well, 
therefore we can determine the densities of the corresponding hyp-hor coverings using the above computation method. 
The density function depends on free parameters $u$ and $p$, and analysing this function we get the minimal density  in 
case $p\approx6.459617$ with $\delta\approx1.268853$. This hyp-hor covering is just locally optimal, because the corresponding 
tiling can not be extended to $\mathbb{H}^3$.
\end{rmrk}
\subsection{Thinnest covering, if the intersection point lies on $A_2P_2$ edge}
In this case, $A_2P_2$ passes through the intersection curve of the balls (see Fig.~5.a), so we can parametrize the intersection 
point of the curve and the edge: $R(r)=(1,x,y,rz_2),r\in [0,1]$. The further computations of this case is very similar to the above two cases. 
We can determine the coordinates of $Q_i,T_i,S_i$ $(i=1,2)$ points, see that the horo- and hyperball cover the edges, 
so the hyp-hor covering is realized, and compute the density of the covering by (4.7). 
The results for tiling $\mathcal{T}_p^{(4,4)}$ (which provides the smallest density in this case) are summarized in the next table. 
\\ \begin{center}
\begin{tabular}{|c|c|c|}
\hline
\multicolumn{1}{|c}{Type of tiling} & \multicolumn{1}{|c|}{$\delta_{min}$} & \multicolumn{1}{|c|}{$r$}\\
\hline
$\mathcal{T}_5^{(4,4)}$      &$1.8383911$ &$0.8114832$ \\
\hline
$\mathcal{T}_6^{(4,4)}$      &$2.3821677$ &$0.7332720$ \\
\hline
$\mathcal{T}_7^{(4,4)}$      &$3.0569894$ &$0.7025236$\\
\hline
\end{tabular}
\end{center}

Finally, summarizing the results so far, we get the following theorems
\begin{theorem}
In $\mathbb{H}^3$, among the hyp-hor coverings generated by simple truncated orthoschemes, the $\mathcal{C}_7^{(3,6)}$ covering configuration 
(see Subsection 4.3) provides the lowest covering density $\approx 1.27297$. The above density is smaller than the so far known lowest covering density 
$\approx1.280$ in the $3$-dimensional hyperbolic space, which was described by L.~Fejes T\'oth and K.~B\"or\"oczky.
\end{theorem}
\begin{theorem}
In hyperbolic space $\mathbb{H}^3$ the function $\delta(\mathcal{C}_p^{(3,6)})$ $(6<p<7,p\in\mathbb{R})$ attains its mimimum in 
case $p \approx 6.459617$, with density $\delta\approx1.268853$, but the corresponding hyp-hor covering can not be extended to 
the entirety of hyperbolic space $\mathbb{H}^3$.
\end{theorem}
We note here, that the discussion of the densest horoball packings in the $n$-dimensional hyperbolic space $n \ge 3$ with horoballs
of different types and hyperballs has not been settled yet.

Optimal sphere packings in other homogeneous Thurston geometries represent
another huge class of open mathematical problems. For these non-Euclidean geometries
only very few results are known (e.g. \cite{Sz07-2}, \cite{Sz14-1}).
Detailed studies are the objective of ongoing research.
%

\noindent
\footnotesize{Budapest University of Technology and Economics Institute of Mathematics, \\
Department of Geometry, \\
H-1521 Budapest, Hungary. \\
E-mail:epermiklos@gmail.com ~ szirmai@math.bme.hu \\
http://www.math.bme.hu/ $^\sim$szirmai}

\end{document}